\documentclass[12pt,notitlepage,twoside,a4paper]{amsart}
 \usepackage{amsfonts}

\usepackage{amsmath,amssymb,enumerate}

\usepackage[usenames, dvipsnames]{color}

\usepackage{xcolor,enumitem}
\usepackage{lipsum}
\usepackage{amssymb}
\usepackage{amsmath,amsthm,comment}
\usepackage{latexsym}
\usepackage{amscd}
\usepackage{psfrag}
\usepackage{graphicx}
\usepackage[latin1]{inputenc}
\usepackage[all]{xy}
\usepackage[mathcal]{eucal}

\everymath{\displaystyle}

\definecolor{NoteColor}{rgb}{1,0,0}

\renewcommand{\textsc}{\textcolor{red}}

\makeatletter
\newcommand\Max{\@tempcnta=\mathcode`\m\relax
\mathcode`\m=\mathcode`\M\max\mathcode`\m=\@tempcnta\relax}

\newcommand\Min{\@tempcnta=\mathcode`\m\relax
\mathcode`\m=\mathcode`\M\min\mathcode`\m=\@tempcnta\relax}
\makeatother

\newtheorem{theorem}{\rm\bf Theorem}[section]
\newtheorem{proposition}[theorem]{\rm\bf Proposition}
\newtheorem{lemma}[theorem]{\rm\bf Lemma}

\newtheorem*{theorem 1}{\rm\bf Proposition 1}
\newtheorem*{theorem 2}{\rm\bf Proposition 2}

\newcommand{\C}{\mathfrak{C}}
\newcommand{\G}{\mathfrak{G}}

\newcommand{\M}{\mathfrak{M}}

\newcommand{\Gx}{\mathfrak{G}_x}

\theoremstyle{definition}

\theoremstyle{remark}

\newtheorem{remarks}[theorem]{\rm\bf Remarks}

\def\interieur#1{\mathord{\mathop{\kern 0pt #1}\limits^\circ}}

\title[Commentary]{A Commentary on Teichm\"uller's paper \\
\emph{Untersuchungen \"uber konforme und quasikonforme Abbildungen} (Investigations on conformal and quasiconformal mappings)}

\author[V. Alberge, M. Brakalova and A. Papadopoulos]{Vincent Alberge, Melkana Brakalova
\\
 and Athanase Papadopoulos}

\thanks{}
   
\date{\today}

\makeindex
\begin{document}

  \maketitle
  \begin{abstract}
  
 This is a commentary on Teichm\"uller's paper \emph{Untersuchungen \"uber konforme und quasikonforme Abbildungen} (Investigations on conformal and quasiconformal mappings) published in 1938. The paper contains fundamental results in conformal geometry, in particular a lemma, known as the Modulsatz, which insures the almost circularity of certain loci defined as complementary components of simply connected regions  in the Riemann sphere, and another lemma, which we call the Main Lemma, which insures the circularity near infinity of the image of circles by a quasiconformal map. The two results find wide applications in value distribution theory, where they allow the efficient use of moduli of doubly connected domains and of quasiconformal mappings. Teichm\"uller's paper also contains a thorough development of the theory of conformal invariants of doubly connected domains.
 
 The final version of this paper will appear in Vol. VII of the \emph{Handbook of Teichm\"uller theory} (European Mathematical Society Publishing House, 2020).

\end{abstract}



\section{Introduction}

The paper \emph{Untersuchungen \"uber konforme und quasikonforme Abbildungen} (Investigations of conformal and quasiconformal mappings)  \cite{T13} is Teichm\"uller's Habilitationsschrift, which he presented in Berlin in 1938, under the supervision of Ludwig\index{Bieberbach, Ludwig (1886--1982)}  Bieberbach.\footnote{Teichm\"uller had  obtained his doctorate in 1935 in G\"ottingen under the supervision of Helmut Hasse.\index{Hasse, Helmut (1898--1979)} The subject was linear operators on Hilbert spaces over the quaternions \cite{T1}. In the meanwhile, he shifted his topics of interest. For details about Teichm\"uller's life, we refer the reader to \cite{abikoff} and \cite{schappacher&scholtz}.} In 1936, Teichm\"uller attended lectures on function theory by Rolf Nevanlinna who was teaching at the University of G\"ottingen as a visiting professor for the academic year. These lectures had a major influence on him. Teichm\"uller remained interested in this topic until the end of his (short) life. His Habilitationsschrift is one of his important papers on complex analysis. He wrote several other papers related to this subject, see \cite{T8, T9, T15, T17, T33}.

The paper  we are commenting on here  is difficult to read. Like several other of Teichm\"uller's paper, while the mathematical ideas it contains are very dense, and the proofs are concise, sometimes  sketchy, Cazacu observes in \cite{Caz}: ``Teichm\"uller's paper(s) are written in a warm, direct style. He emphasizes aims, ideas, difficulties, clearly explains notations, methods, and compares his results with others concerning similar problems.'' All throughout the paper, Teichm\"uller gives credit to  the works of other mathematicians which he uses, improves or has been motivated by. Another characteristic of this  paper is that  its spirit and the proofs it contains  are purely geometric. We may quote here a passage from the paper, during a computation he makes of the module function $\Phi(P)$. Teichm\"uller writes  (\S 2.1): 
\begin{quote}\small
[\ldots] Explicitly calculating this requires the introduction of elliptic functions and one can observe a connection between the function $\Phi(P)$ and the values of the {\it elliptic modular function} for purely imaginary periodic ratios. In the interest of preserving the purity of the method used and to avoid tedious calculations, we will not make use of this connection at all; rather, we derive everything we wish to know about the function $\Phi(P)$ from its geometric definition.
   \end{quote}

The paper has seven sections. The first six are concerned with the conformal geometry of doubly- and simply-connected domains (the latter with marked points),  extremal problems for conformal invariants,\index{conformal invariant}\index{conformal invariant!extremal problem} and quasiconformal mappings.\index{quasiconformal mapping} Teichm\"uller studies  properties of modules of quadrilaterals and ring domains and reduced modules\index{reduced module} of simply-connected domains. He often relies on applications of the  length-area method. One of the main motivations for his investigations is to give a geometric proof of what he calls the \emph{Main Lemma} (see Item (\ref{2}) and Theorem \ref{ML} below), stated at the beginning of the paper, of which an improved version is now known as the Teichm\"uller--Wittich--Belinski\u{\i} theorem.\index{Teichm\"uller--Wittich--Belinski\u{\i} theorem}\index{Theorem!Teichm\"uller--Wittich--Belinski\u{\i}} This result  was motivated by the developments of (at that time) new value distribution theory\index{value distribution}. Indeed, Teichm\"uller says that he reached this statement while we was working with Hans Wittich on a class of Riemann surfaces introduced by Egon Ullrich. (Such surfaces are now called surfaces with finitely \emph{periodic ends}). Teichm\"uller also writes that he first tried to prove this lemma without success by using analytical methods developed  by Ahlfors in \cite{Ahlfors:thesis}  and \cite{Ahlfors-Acta}. Throughout Teichm\"uller's paper, there is no explanation about how this result could be used and we have to quote another paper by Teichm\"uller. In \cite{T20}, Teichm\"uller writes: 
\begin{quote}
In order to study a conformal mapping from a Riemann surface $\mathfrak{A}$ to another Riemann surface $\mathfrak{B}$, one constructs an explicitly known quasiconformal mapping from $\mathfrak{A}$ to some Riemann surface $\mathfrak{C}$ [...] so that $\mathfrak{B}$ and $\mathfrak{C}$ are closely related and so that the best possible estimate from above for the dilatation quotients viewed as a function of points of $\mathfrak{C}$ is known. [...] One then determines the properties attached to the mappings from $\mathfrak{C}$ to $\mathfrak{B}$ whose dilatation quotients satisfy this estimate. Thanks to the known mapping $\mathfrak{A}\rightarrow  \mathfrak{C}$, each of these properties gives a statement about the mapping $\mathfrak{A}\rightarrow  \mathfrak{B}$ to be studied.
\end{quote}

This method and therefore the Main Lemma has been first applied by L\^{e} V\u{a}n  \cite{Le} to solve a particular case of the so-called \emph{Nevanlinna inverse problem}. We refer to the paper \cite{P-value} in the present volume and to the paper \cite{DGP} for an exposition of this problem, to which, later on, Drasin provided a full solution in \cite{Drasin1977}; see also the monograph \cite{GolOst}. The last section of Teichm\"uller's paper \cite{T13} is a  contribution to the  type problem,\index{type problem} in which the author gives a negative answer to a conjecture made by  Nevanlinna. The methods rely to some extent on the quasiconformal techniques developed in the previous sections, but are not based on the Main Lemma.

Let us add that Teichm\"uller's main interest in the paper we are concerned with here is the question of uniformization and conformal representation. He  writes (\S 5.1): ``
[\ldots] it seems necessary to stress that I do not investigate the concept of  module of a ring domain for its own sake, but that the focus is on an investigation of the problem of conformal mappings and uniformization."

Let us give a list of other main results and tools that are introduced in this paper.
 \begin{enumerate}

  \item A study of the conformal geometry of annuli and quadrilaterals: monotonicity of modules, regions that are solutions of extremal problems, the introduction of new conformal objects: conformal radius, reduced module,\index{reduced module} reduced logarithmic area,\index{reduced logarithmic area} and others.

   \item  \label{4} A distortion theorem which is an improvement of the main distortion lemma in Ahlfors' thesis, obtained (unlike Ahlfors' proof) without using differential equations but using only geometric methods. 
                
    \item An important lemma known today as Teichm\"uller's Modulsatz\index{Modulsatz} which is used in his proof of the Main Lemma (see (\ref{2}) or \ref{ML}). 
                The lemma gives a sufficient condition for a set that lies between two doubly-connected domains in the complex plane to be close to a circle.                     

 \item \label{1} The use of a class of quasiconformal mappings\index{quasiconformal mapping} $w=w(z)$ of the complex plane which is more general than the class used today: the dilatation quotient of such a map at each point (that is, the quotient of the great to the small axis of the infinitesimal ellipse whic is the image of an infinitesimal circle) is not required to be bounded by a uniform constant, but it is assumed that it is bounded by some specific function of $\vert z\vert$ with controlled growth.      
 
  \item  \label{2} An asymptotic result for quasiconformal mappings, namely, a bound as in (\ref{1}), that insures that $\vert w\vert \sim \mathrm{const.} \vert z\vert$ as $z\to\infty$ (a property we call ``circularity at infinity" of the image of circles by such quasiconformal mappings). Following Teichm\"uller, we shall refer to this result as the Main Lemma.\index{Main Lemma!Teichm\"uller}\index{Teichm\"uller Main Lemma}  

     \item A contribution to the type problem,\index{type problem} and more precisely, a condition for a Riemann surface to be of hyperbolic type. Teichm\"uller uses for this purpose the notion of line complex, a combinatorial device to encode a Riemann surface that is a branched cover of the sphere. At the same time, he settles (by the negative) a question posed by Nevanlinna on the type problem.

             \end{enumerate}   

We have divided this commentary into sections which highlight the main results we stated above. More precisely, in Section \ref{s:mod} we will recall the geometric tools used by Teichm\"uller, such as the module of a ring domain, the reduced module of a simply-connected domain and the module of a quadrilateral. In Section \ref{s:ahlfors}, we will comment on  the improvement made by Teichm\"uller to the Ahlfors distortion theorem. In Section \ref{s:Modulsatz}, we will present the main ideas of the proof of the Modulsatz theorem along with some applications.  In the last section, Section \ref{s:main}, we will recall the notion of quasiconformal mappings that Teichm\"uller uses and we will explain the Main Lemma and its improvements.  Let us point out that  we have decided to comment on the last section of \cite{T13}, that is the contribution to the type problem, in another chapter of the present volume, \cite{ABP-type}.

  
  \section{Module theorems}\label{s:mod} 
  
In the first two sections of his paper, Teichm\"uller develops from first principles  properties of  modules of  ring domains using logarithmic area and logarithmic length and modules of quadrilaterals, and he defines and studies  properties of reduced modules\index{reduced module} of  simply-connected domains. The references to his paper for what concerns reduced module, in the monographs by Wittich \cite{Wit68} and  Jenkins \cite{JJ58},  suggest that  he was the first to develop this kind of results. We can also mention the book by K\"unzi \cite{kunzi} which contains the main results of  Teichm\"uller  from \cite{T13} that we are commenting on here. In his study of moduli of ring domains, Teichm\"uller, as he mentions in the introduction, was motivated by an estimate of Ahlfors, which he uses after applying a suitable auxiliary conformal map. Since this involves only considerations of the point $\infty$, he was led  to study conformal mappings of neighborhoods of this point, excluding this point itself, that is, of annuli.

Let us recall in detail the notions introduced and  used in  Teichm\"uller's paper  \cite{T13}.

Teichm\"uller starts in \S1.1 by recalling the definition of a \emph{ring domain}, that is, a doubly-connected domain bounded by two continua. 

In \S1.2, he shows that any ring domain $\G$ can be mapped conformally onto an annulus $r<|w|<R,$ where $0 < r<R< \infty$. The proof uses the uniformization theorem (Teichm\"uller refers to Koebe).  He then defines  the module\index{module!of a doubly-connected domain} $M$ of $\G$  as $\log \frac{R}{r}$. This is a conformal invariant\index{conformal invariant} and the definition is equivalent to a well-known classical definition based on extremal length, see \cite{Ahlfors-CI, Ahlfors1964}.  Let us point out that Teichm\"uller's definition of module differs by a multiplicative ratio of $2\pi$ from the definition of the module that involves the notion of extremal length. We also note that conformal invariants of ring domains were already studied by Schottky. The latter proved in his doctoral dissertation \cite{Schottky} (1877) that two annuli $\{r_1<\vert z\vert <r_2\}$ and $\{r'_1<\vert z\vert <r'_2\}$ are conformally equivalent if and only if we have $\displaystyle \frac{r_2}{r_1}=\displaystyle \frac{r'_2}{r'_1}$.  

In \S 1.3,  Teichm\"uller introduces the  \emph{logarithmic length} of a curve $\gamma$ in $\mathbb{C}\setminus\left\lbrace 0\right\rbrace$,  defined as \[\int\limits_{\gamma}|d\log z|=\int\limits_{\gamma} \dfrac{|dz|}{|z|^2},\]
 and the \emph{logarithmic area} $F$ of a domain $\G \subset \mathbb C \setminus \left\lbrace 0\right\rbrace $,  defined as \[\iint\limits_G \framebox{$d\log z$}= \iint\limits_G \dfrac{\framebox{$d z$}}{|z|^2},\]
where, ${\framebox{$d z$}}=dxdy$.  These two quantities  represent the Euclidean length and area of the  transformed, to the $\log z$-plane, original curve and domain, where the $\log z$-plane is the rectangle in the Euclidean plane which is the image of the annulus $r<|w|<R$ by the logarithm map, after the annulus has been cut along the interval $[-R,-r]$. 
If  $M$ and $F$ denote respectively the module and the logarithmic area of $\G$, Teichm\"uller shows that
\begin{equation}\label{eq:comparison-log}
2\pi M\leq F .
\end{equation}  Let us say a few words about Teichm\"uller's proof of (\ref{eq:comparison-log}). This proof uses the so-called \emph{length-area method}, a method that he refers to in the introduction as the ``Gr\"otzsch's and Ahlfors' methods''  Let us recall that the length-area method involves comparison of Euclidean lengths of curves and  area in the logarithmic plane $w=\log z$-plane,  followed by  a repeated integration and use of the Cauchy--Schwarz inequality.  The concept of extremal length\index{extremal length} takes root in the use of this method which was initiated by Beurling in his thesis \cite{Beurling-thesis} (see also the end of this section), and further developed by Beurling and Ahlfors\footnote{Cf. the historical notes in Ahlfors' book \emph{Conformal invariants} \cite[p. 81]{Ahlfors-CI}.} in the 1940s by working not only  with the Euclidean metric but with a whole class of conformal metrics 
(see \cite{Ahlfors1946} and \cite[p. 50]{Ahlfors-CI}). According to Ahlfors,\footnote{Cf. L. V. Ahlfors', Collected works \cite{Ahlfors-Collected}, Vol. 1, p. 1} his first used of the length-area  method  appears in his thesis (see also \cite{Ahlfors-CRAS,Ahlfors-Lindelof,Ahlfors:thesis})  and was inspired by the book \emph{Funktionentheorie} (1922) by Hurwitz and Courant \cite{HC}.   

In \S1.4 and \S1.5, Teichm\"uller presents an elegant proof  of the following two  properties of the module, and, at the same time, a solution of two extremal problems:
\begin{proposition}\label{comparison} If a ring domain $\G^\prime$ with module $M^\prime$ is a subset of a ring domain $\G$ with module $M$ such that $\G^\prime$ separates the two complementary continua of $\G$, then $ M'\leq M.$ {\it Equality holds if and only if $\G'=\G.$} \end{proposition}
\begin{proposition} \label{subadd} If a ring domain $\G$ contains two disjoint ring domains $\G'$ and $\G''$ each of which separates the complementary continua of $\G$ and has modules $M'$ and $M''$  respectively, then $ M'+M''\leq M$. Equality holds if and only if $\G=\{r<|z|<R\},\ \G'=\{r<|z|<\rho\}$ and $\G''=\{\rho<|z|<R\}.$ \end{proposition} 
This kind of module inequalities and of solutions of extremal problems involving ring domains were first studied by Gr\"otzsch in \cite{Groetzsch1928}. The reader may also refer to the chapter \cite{AP1} in the present volume containing an exposition of these results.

In \S 1.6, Teichm\"uller starts by recalling the notion of \emph{conformal radius}\index{conformal radius}, which he calls ``mapping radius.'' We recall that the conformal radius of a simply-connected domain $\G$ of $\mathbb{C}$ containing $0$ is the positive number $R$ for which $\G$ is mapped conformally  onto the unit disc such that $0$ is mapped onto $0$ and the derivative at $0$ is equal to $\frac{1}{R}$. He continues by defining the notion of \emph{reduced module} of a simply-connected domain $\G$ containing $z=0$ as follows. Let $\rho >0$ be small enough in order to  have the disc of radius $\rho$ centered at $z=0$  contained in $\G$ and let $M_{\rho}$ be the module of the ring domain obtained from $\G$ by removing such a disc.  By means of what Teichm\"uller calls a ``distortion theorem'' (a term which refers to one of Koebe's theorems), he obtains the following estimate:
\begin{equation}\label{redmodkoebe}
|M_{\rho}+\log \rho-\log R|\leq \dfrac{2\rho}{R-4\rho}.
\end{equation}
Then, the reduced module $\widetilde{M}$ of $\G$ at $z=0$ is defined as  $\lim\limits_{\rho\to 0}(M_\rho+\log \rho)$ and by means of (\ref{redmodkoebe})  is equal to  the logarithm of the conformal radius of $\G$. The definition is then extended to the reduced module\index{reduced module} at any point of $\G$, and by the transformation $z\mapsto\dfrac 1z$, if  $\infty\in \G$, it is also extended  at $z=\infty$. In  the special case where  $\G'=\{0<|z|<R\}$ and  $\G''=\{ |z|>R\}$,  we have $\widetilde{M}' =\log R$ and $\widetilde{M}''=-\log R,$ respectively,  thus $\widetilde{M}'+\widetilde{M}''=0$.  

In \S1.7, Teichm\"uller defines the \emph{reduced logarithmic area} of a simply-connected domain $\G$ containing $z=0$ but not $z=\infty$ by the formula
 $$\widetilde{F}=F_{\rho}+ 2\pi \log \rho,$$ where $\rho$ is small enough so that the disc $D_{\rho}$ of radius $\rho$ centered at $z=0$ is contained in $\G$, and $F_{\rho}$ denotes the logarithmic area of the ring domain $\G\setminus D_{\rho}$. He then obtains an inequality similar to (\ref{eq:comparison-log}), namely, he proves that for a simply-connected domain $\G$ containing $z=0$ of reduced module $\widetilde{M}$ and of reduced logarithmic area $\widetilde{F}$, we have 
\begin{equation}\label{eq:comparison-log-reduced}
2\pi \widetilde{M}\leq \widetilde{F},
\end{equation}
and this equality holds if and only if $\G=\left\lbrace \left| z \right| < e^{\widetilde{M}} \right\rbrace$. 

Teichm\"uller's proof involves the length-area method but he points out that it can be done by what he calls   ``Bieberbach's method.'' Teichm\"uller justifies such a preference by saying that he only wants to apply geometric methods.

In \S1.8, as a conclusion of its first section, Teichm\"uller deduces  two results on simply-connected domains which are analoguous to Propositions \ref{comparison} and \ref{subadd}. He also obtains the following result which is of special importance since it is used in the proof of the special Modulsatz, see \S 4 below.
\begin{proposition}\label{lemma:special-modulsatz}
Let $\G'$ and $\G''$ be two disjoint simply-connected domains with $\G'$ containing $z=0$ and $\G''$ containing $z=\infty$ and having reduced modules $\widetilde{M}'$ and $\widetilde{M}''$ respectively. Then
$$
\widetilde{M}'+\widetilde{M}''\leq 0.
$$
Equality holds only when $\G'=\left\lbrace |z|<e^{\widetilde{M}'} \right\rbrace$ and $\G''=\left\lbrace |z|>e^{\widetilde{M}'} \right\rbrace$.
\end{proposition}

We now review extremal domains studied by Teichm\"uller in \S 3 and \S 4.
 
  In \S 2.1--\S2.3, Teichmüller 
 works with a particular class of ring domains, namely,  for $P>1$, the class of ring domains that separate the unit disc from $z=\infty$ such that the distance of the unbounded component of their corresponding complementary region from the origin is at most $P$.  
  
Let  ${\G}_P$ be  the exterior of the unit disc cut along the real axis from $z=P$ to $z=\infty$.  Teichm\"uller proves several properties of ${\G}_P$, including the following:
\begin{proposition}\label{proposition-grotzsch-module} Let  $\log\Phi(P)$ be the module of ${\G}_P$. 
Then,

 (1)  $\Phi(P) > P$.

(2)  $\Phi(P)$ is a continuous strictly monotone increasing function of $P$, and $\Phi(P)$ takes all values between $1$ and $\infty$.

 (3) $\Phi(P)/ 4P\to 1$  in an increasing manner as $P\to \infty$. 

 (4) Suppose that for $P>1$ the ring domain $\G$ belongs to the above class. Let $M$ be the module of this domain. Then,
$$ M\leq \log \Phi(P)$$
and equality holds if and only if $\G$ is equal to some  ${\G}_P$ up to a rotation about the origin.

\end{proposition}
Let us make a few remarks on these results. 

\begin{remarks}
 1.-- Items (1), (2) and (3) are consequences of  Inequality (\ref{redmodkoebe}) above, with $R=4$. 
 
 2.-- The value $\log 4$ is the reduced module of the whole plane cut along the real axis from $1$ to $\infty$.
 
 3.--  Teichm\"uller's proofs of these results do not involve use of elliptic integrals and rely only on geometric methods. The reason, as he says, is that he wanted to preserve ``purity'' and ``avoid tedious calculations.'' 
  
  4.--  Proofs of such results by means of elliptic integrals can be found in Hersch's thesis \cite{herschthesis} and also in the book by Lehto and Virtanen \cite{LV73}. 
  
  5.--  The result in item (4) is referred to in  \cite{LV73} as the \emph{Gr\"otzsch module theorem}---the\index{Gr\"otzsch module theorem}\index{module theorem!Gr\"otzsch}\index{Theorem!Gr\"otzsch module} extremal domain  ${\G}_P$ is called the \emph{Gr\"otzsch extremal domain}\index{Gr\"otzsch extremal domain}\index{extremal domain!Gr\"otzsch}---and Teichm\"uller proves it in two different geometric ways. He also explains how one of them can be used to prove ``Koebe's theorem with the exact constant $\frac{1}{4}$.'' 
  
 Teichm\"uller notes at the beginning of \S 2.1 that  this extremal problem has been already solved by Gr\"otzsch in \cite{Groetzsch1928}, using the so-called \emph{Gr\"otzsch method of strips}.   A proof based on the \emph{extremal length}\index{extremal length} is contained in \cite{Ahlfors-QC} and in \cite{LV73}. 
   
 6.-- Teichm\"uller uses these results in \cite{T31} in order to prove that ``a quasiconformal mapping of the unit circle $\left| z \right| < 1$ onto itself, where all boundary points remain fixed, shifts the center $z=0$ by at most $2\left( \left(\sup D(z)\right)^{\frac{1}{2}}-1\right)$.'' We refer to \S \ref{s:main} below for the definition of a quasiconformal mapping and the definition of the associated dilatation quotient $D(z)$.
 
 \end{remarks}


%
%
In \S 2.4, 
Teichm\"uller solves another extremal problem, namely, he proves that for $1<P_1<P_2$, the annulus $1<\vert z\vert <P_2$ with a slit along the segment $\left[P_1,P_2\right]$ has the largest module amongst all ring domains separating their $\left| z \right| =1$ boundary component from $\infty$ and such that the outer boundary component lies in the annulus $r<\vert z\vert <P_2$, for some $1<r\leq P_1$. After proving this extremal property, Teichm\"uller concludes this subsection by a result which is a weak version of the Modulsatz (see Theorem \ref{Modulsatz} below for the statement). 

In \S 2.5, Teichm\"uller leaves the investigation of extremal problems in order to study normal families  between doubly-connected domains, where the base domain is the annulus $1<\vert z\vert <R$.

He then returns in \S 2.6 (and until the end of \S 2) to the study of another extremal problem and proves what is often referred to   as   \emph{Teichm\"uller's module theorem}.\index{Teichm\"uller module theorem}\index{module theorem!Teichm\"uller}\index{Theorem!Teichm\"uller module} (This is a terminology used in Lehto and Virtanen's monograph  \cite{LV73}.)  We recall that this result says that  among all doubly-connected domains  that separate the points $0$ and $\rho e^{i\varphi}$ from $Pe^{i\theta}$ and $\infty$, the one with the largest module is (up to a rotation) the so-called \emph{Teichm\"uller extremal domain},\index{extremal domain!Teichm\"uller} that is,  the complex plane cut along the real axis from $-\rho$ to $0$ and from $P$ to $\infty$. Teichm\"uller denotes the module of the latter by $\log \Psi \left( \frac{P}{\rho} \right)$. The proof of this property is sketchy (it is not as detailed as the one for the corresponding result on Gr\"otzsch's module). A more detailed proof can be found for instance in  \cite{LV73}, \cite{Ahlfors-QC} and  \cite{Ahlfors-CI}. 

Teichm\"uller then compares the functions $\Psi$ and $\Phi$. He obtains, using conformal mappings, 
$$
\Psi \left(\frac{P}{\rho} \right) = \Phi \left( 1+2\frac{P}{\rho}\left( 1+\sqrt{1+\frac{P}{\rho}} \right)\right)
$$
and
$$
\Psi \left(\frac{P}{\rho} \right) = \Phi \left( \sqrt{1+\frac{P}{\rho}} \right)^2.
$$
From these two relations and Item (3) of Proposition \ref{proposition-grotzsch-module} he deduces the two estimates 
\begin{equation}\label{estim te1}
\Psi \left(\frac{P}{\rho} \right) \sim 16\frac{P}{\rho} \textrm{ as } \frac{P}{\rho}\rightarrow \infty
\end{equation}
and
\begin{equation}\label{estimate2}
\Psi \left( \frac{P}{\rho}\right) < 16\left( \frac{P}{\rho}+1 \right).
\end{equation}

In the last subsection, \S 2.8,  as a consequence of the extremal problem, Teichm\"uller deduces that any doubly-connected domain  of the plane separating $0$ from $\infty$ whose module is at least equal to $e^{\pi}$ always contains a circle of center $0$.

We conclude our comments on \S 2 by adding that, as Teichm\"uller writes, the study of extremal domain and its module are related to Ahlfors' proof of  Denjoy's conjecture. As we shall see, this relation leads to an improvement of the  distortion lemma of Ahlfors. 

  Teichm\"uller starts \S 3 by recalling different conformal invariants associated with  a quadrilateral, namely, the \emph{cross-ratio}, the \emph{harmonic measure}\index{harmonic measure} and the \emph{module}. He concludes \S 3.1 by explaining how these notions are related to each other.
 Let us recall here that the module\index{module!of a quadrilateral} of a quadrilateral is the ratio of the sides of a Euclidean rectangle that is conformally equivalent to it.  
 
  Teichm\"uller gives in \S 3.2 an estimate of this module that he qualifies as ``simple''  and which is at the basis of the length-area method.\index{length-area method} We give below the statement, following Teichm\"uller's notation.
%
\begin{proposition}\label{prop-extremal-distance}
Let $\mathfrak{V}$  be a quadrilateral with sides $\mathfrak{a},\mathfrak{b},\mathfrak{c},\mathfrak{d}$ and let it be mapped conformally onto the rectangle $0<u<a, \, 0<v<b$ in the $w=u+iv$-plane, with $\mathfrak{a}, \mathfrak{b},  \mathfrak{c}, \mathfrak{d}$ mapped respectively to $v=0, \ u=a,\ v=a,\ u=0$. Let $\beta$ be the infimum of the lengths of all curves in $\mathfrak{V}$ joining $\mathfrak{a}$ and $\mathfrak{c}$, and let {\bf $F$} be the area of $\mathfrak{V}$. Then
$$
\dfrac{a}{b}\leq\dfrac{F}{\beta^2}.
$$
Equality holds if and only if $\mathfrak{V}$ is a rectangle with sides $\mathfrak{a},\mathfrak{b},\mathfrak{c},\mathfrak{d}$.
\end{proposition}

 In \S 3.2, Teichm\"uller derives some consequences from this result. For instance, using the same notation as above, he proves that 
 \begin{equation}\label{inequalityThersch1}
 \alpha \beta \leq F
 \end{equation}
 with equality only for rectangles. He also proves that such an inequality is also satisfied for  quadrilaterals on a Riemannian surface. 
 He then observes that 
 \begin{equation}\label{inequalityThersch2}
 \min\left\lbrace \alpha , \beta \right\rbrace^2 \leq F,
 \end{equation}
 with equality only for squares.\footnote{These two results were proved later by Besicovitch in  \cite{Bes52} in a manner which according to Jenkins in \cite{JJ58} is a ``length-area proof of primitive type.'' Besicovitch in his paper does not mention Teichm\"uller though he refers to L\"owner. Let us add that using the same strategy as Teichm\"uller, it is not difficult to prove that one obtains the same inequalities as (\ref{inequalityThersch1}) and (\ref{inequalityThersch2}) by replacing the Euclidean lengths $\alpha$ and $\beta$ by appropriate extremal lengths. The inequalities so obtained are attributed by Hersch in \cite{hersch55} and Jenkins in \cite{JJ58} to Teichm\"uller.}
%
  
Teichm\"uller concludes this subsection by generalizing the notion of extremal distance\index{exremal distance} introduced by Beurling in his thesis, \cite{Beurling-thesis}, for a pair of disjoint boundary arcs, and proves by means  of Proposition  \ref{prop-extremal-distance} that in the case of a quadrilateral, the extremal length coincides with the reciprocal of the associated module. 

One has to note here that Teichm\"uller's notion of  extremal distance is different from the notion of extremal distance of a pair of disjoint boundary curves introduced later on by Ahlfors and Beurling.\footnote{Indeed, it is not difficult to see that  Teichm\"uller only deals with conformal metrics that are pullback of the Euclidean metric by conformal mappings.} 
  
We are now ready to talk about Ahlfors' distortion result and its improvement by Teichm\"uller.

\section{Ahlfors' distortion theorem}\label{s:ahlfors}

 Ahlfors proves the Denjoy conjecture in his doctoral dissertation, written under the supervision of Nevanlinna and published in \cite{Ahlfors:thesis}.\index{Denjoy conjecture}\index{conjecture!Denjoy} The result is announced in \cite{Ahlfors-CRAS}, and it also appears in \cite{Ahlfors-Lindelof}.  This conjecture says that the number of finite asymptotic values of an entire function of order $k$ is at most $2k$. Making this statement precise needs the definition of an appropriate notion of ``order" and ``asymptotic value" of an entire function, and this is done in the setting of Nevanlinna's theory (or value distribution theory);\index{value distribution} see the paper \cite{P-value} in the present volume. As mentionned earlier in the paper,  Ahlfors' approach to the Denjoy conjecture was completely new and based on the 
length-area method.  From such a method he obtains two ``main inequalities," one of which is  called the ``Erste Hauptungleichung'' (first main inequality) and is the key result for the proof of the Denjoy conjecture. It is now known as the Ahlfors distortion theorem.\index{Ahlfors distortion theorem}\index{distortion theorem!Ahlfors}  
 Regarding this result, Ahlfors writes in \cite[Vol. 1, p.1]{Ahlfors-Collected}:
\begin{quote} [...]  In my thesis \cite{Ahlfors:thesis} the lemma on conformal mapping has become the main theorem in the form of a strong and explicit inequality or distortion theorem for the conformal mapping from a general strip domain to a parallel strip, together with a weaker inequality in the opposite direction. [...] A more precise form of the first inequality was later given by O. Teichm\"uller.
\end{quote}
One can also point out that Jenkins and Oikawa in \cite{JO71} gave another proof of the two Ahlfors inequalities by means of the extremal length method. Their proof of the first one was described by Ahlfors as a ``virtually trivial proof.''

Following Teichm\"uller's exposition (from \S 3.5 to \S 3.7), we now present  Ahlfors' distortion theorem along with the main ideas of its proof given by Teichm\"uller. For this purpose, we first need to set the framework.  

 Let $\G$ be a strip domain, that is, a simply-connected domain of the $z=(x,y)$-plane with  two marked accessible  boundary points (that is, points that can be joined by a continuous curve to an interior point) $\mathfrak{r}_1$ and $\mathfrak{r}_2$ satisfying $\Re \mathfrak{r}_1< \Re \mathfrak{r}_2$.

   For every $x$ in the open interval $\Re \mathfrak{r}_1<x< \Re \mathfrak{r}_2,$ the strip domain $\G$ is divided by the line $\Re z = x $  into at most countably many subsets of which one, denoted by $\Gx$, contains $\mathfrak{r}_1$ on its boundary. This set $\Gx$ has also the property that the last piece of every curve in $\G$ ending at $\mathfrak{r}_1$ lies in $\Gx$. The interior points of $\G$ which are boundary points of $\Gx$ satisfy $\Re z = x$ and are divided into countably many cross-cuts $\mathfrak{S}_i$ of $\G$, that is, simple arcs whose endpoints lie in the boundary of $\G$ and whose interior points lie in the interior of $\G$. Each of these cross-cuts $\mathfrak{S}_i$ partitions $\G$ into two pieces, one containing $\Gx$ and another one which we will call for now $\G_i$. The strip $\G$ is then the disjoint union of $\Gx$, the $\mathfrak{S}_i$, and $\G_i$. Note that all cross-cuts $\Re z = x$ of $\G$ that do not lie on the boundary of $\Gx$ are erased, so that  $\G_i$ may contain points with $\Re z < x$. Whichever $\G_i$ contains $\mathfrak{r}_2$ as boundary point, and also contains the end piece of each curve in $\G$ terminating in $\mathfrak{r}_2$, will be separated from $\G_x$ by the cross-cut $\mathfrak{S}_i=\mathfrak{S}_x$. 
   
   Let $\Theta(x)$ be the length of $\mathfrak{S}_x$. 
   As Nevanlinna points out  in \cite[p. 93]{N-analytic}, such a function ``[...] is in general of a complicated nature,"  but Teichm\"uller proves that it is lower semi-continuous and therefore measurable. 
%
By the uniformization theorem, one maps the domain $\G$ onto the parallel strip $\left\lbrace w=u+i v \mid 0< v < B\right\rbrace$ such that the boundary points $\mathfrak{r}_1$ and $\mathfrak{r}_2$ are sent respectively to $-\infty$ and $\infty$. For $\Re \mathfrak{r}_1 <x< \Re \mathfrak{r}_2$,   the cross-cut $\mathfrak{S}_x$ is thus mapped onto an arc $L_x$ that connects the two boundary straight lines of the strip.  
 
 We set  
 \[u_1(x) = \min\left\lbrace \Re w \mid w\in L_x \right\rbrace\] and 
 \[u_2 (x)=  \max\left\lbrace \Re w \mid w\in L_x \right\rbrace.\]
 
  Then, as recalled by Teichm\"uller, Ahlfors' distortion theorem\index{Ahlfors distortion theorem}\index{distortion theorem!Ahlfors} says that 
 \begin{equation}\label{eq:ahlfors-distortion}
 \frac{u_1 (x^{\prime\prime})-u_{2}(x^{\prime})}{B} > \int_{x^{\prime}}^{x^{\prime\prime}}{\frac{dx}{\Theta(x)}}-4
 \end{equation}
whenever for $\Re \mathfrak{r}_1 < x^{\prime}<x^{\prime\prime}<\Re \mathfrak{r}_2$, $\int_{x^{\prime}}^{x^{\prime\prime}}{\frac{dx}{\Theta(x)}}>2$. In order to prove this,
 Teichm\"uller first proves by means of the length-area method that
\begin{equation}\label{eq:ahlfors-distortion-teichmueller1/2}
\forall x\in\left( \Re \mathfrak{r}_1 ,  \Re \mathfrak{r}_2\right), \; \int_{x^{\prime}}^{x^{\prime\prime}}{\frac{dx}{\Theta(x)}}\leq \frac{a}{b},
\end{equation}
where $a$ and $b$ are  the lengths of the sides of the rectangle $\mathcal{R}(a,b)$ that is biholomorphic to the quadrilateral $\mathcal{Q}(x^{\prime},x^{\prime\prime})$ bounded by $\mathfrak{G}$ and the cross-cuts $\mathfrak{S}_{x^{\prime}}$ and $\mathfrak{S}_{x^{\prime\prime}}$  and for which the sides of length $a$ are images  of the cross-cuts by such a biholomorphism.  
 Then, he continues by applying suitable conformal mappings as follows. First, he  maps the rectangle $\mathcal{R}(a,b)$ conformally onto the upper half of an annulus of module $\pi\frac{a}{b}$. Then he maps the quadrilateral $\mathcal{Q}(x^{\prime},x^{\prime\prime})$ onto the upper half-part of a symmetric doubly-connected domain that separates the points $0$ and $1$ from $e^{\pi\frac{u_{1}(x^{\prime\prime})-u_{2}(x^{\prime})}{B}}$ and $\infty$. 
 
  From the extremal problem he proves earlier, namely, the so-called Teichm\"uller module theorem,\index{Teichm\"uller module theorem}\index{module theorem!Teichm\"uller}\index{Theorem!Teichm\"uller module} he obtains 
 \begin{equation}\label{eq:ahlfors-distortion-teichmueller2/2}
 \frac{a}{b}\leq \frac{1}{\pi}\log\Psi\left( e^{\pi\frac{u_{1}(x^{\prime\prime})-u_{2}(x^{\prime})}{B}} \right).
 \end{equation}
By using the inverse function of $\frac{1}{\pi}\log\Psi(e^{\pi x})$, inequalities (\ref{estimate2}), (\ref{eq:ahlfors-distortion-teichmueller1/2}) and (\ref{eq:ahlfors-distortion-teichmueller2/2}), Teichm\"uller finally obtains
\begin{equation}\label{eq:ahlfors-distortion-teichmueller}
\frac{u_{1}(x^{\prime\prime})-u_{2}(x^{\prime})}{B}>\int_{x^{\prime}}^{x^{\prime\prime}}{\frac{dx}{\Theta(x)}}-\frac{4\log 2}{\pi}-\frac{1}{\pi}\log\frac{1}{1-8e^{-\pi\int_{x^{\prime}}^{x^{\prime\prime}}{\frac{dx}{\Theta(x)}}}}.
\end{equation}
It now becomes elementary to verify that such an inequality implies Ahlfors' distortion theorem. 


Let us conclude this section by adding that the main steps of Teichm\"uller's proof can be found  in  \cite[p. 97]{N-analytic} and \cite[p. 76]{Ahlfors-CI}. 

\section{The Modulsatz}\label{s:Modulsatz}

We now comment on the Modulsatz  and its consequences. 

At the beginning of \S 4, Teichm\"uller explains that he is interested in estimating the location of the points that lie between two disjoint ring  subdomains 
of a given annulus $\G =\left\lbrace z \mid r< \left| z \right| < R\right\rbrace$. 
 For this purpose, he states what he calls the Modulsatz and which says the following: 
\begin{theorem}[The Modulsatz]\label{Modulsatz}\index{Modulsatz} For every $\varepsilon>0$, there exists $\delta >0$ such that for any two disjoint ring subdomains $\G^{\prime}$ and $\G^{\prime\prime}$ of $\G$ where $\G^{\prime}$ separates $0$ from  $\G^{\prime\prime}$ and with modules $M$ and $M^{\prime}$ respectively, if
$$
M'+M''\geq \log\frac{R}{r}-\delta,
$$
then each point separated from $0$ by $\G'$ and from $\infty$ by $\G''$ belongs to the annulus
\begin{equation}\label{estimate}
\log r+M'-\varepsilon\leq\log |z|\leq \log R-M''+\varepsilon.
\end{equation}
\end{theorem}
As he explains, such a result is a consequence of the following:
\begin{theorem}[Special Modulsatz] \label{Specialmodulsatz}\index{special Modulsatz}\index{Modulsatz!special}  For every $\varepsilon >0$, there exists $\delta=\delta(\varepsilon)$ such that for any two disjoint simply connected domains of the Riemann sphere $\widetilde{\G}^{\prime}$ and $\widetilde{\G}^{\prime\prime}$  containing respectively $0$ and $\infty$ and of reduced modules respectively $\widetilde{M}'$ and $\widetilde{M}^{\prime\prime}$, if 
$$
\widetilde{M}'+\widetilde{M}''\geq -\delta,
$$ 
then   the complement of $\widetilde{\G}^{\prime}\cup\widetilde{\G}^{\prime\prime}$ is contained in the circular ring
$$
\widetilde{M}'-\varepsilon \leq \log |z|\leq -\widetilde{M}''+\varepsilon.
$$
\end{theorem}
Indeed,  one has just to apply Theorem \ref{Specialmodulsatz} to the regions obtained from $\G^{\prime}$ and $\G^{\prime\prime}$ by adjoining to them their corresponding complements containing $0$ and $\infty$ respectively. 

Since the Modulsatz is a consequence of the special Modulsatz,\index{special Modulsatz}\index{Modulsatz!special} Section 4 of Teichm\"uller's paper is mainly  devoted to the proof of the latter.

 Teichm\"uller gives two proofs. The first one is a proof by contradiction which is based on a normal family argument and the Rouch\'{e} theorem but as Teichm\"uller points out it is ``purely existential,'' in the sense that it only assures the existence of $\delta=\delta(\varepsilon)$ without specifying how it depends on $\varepsilon$. His second proof gives the relationship between $\delta$ and $\varepsilon$ and is 
among the major geometric and more sophisticated achievements of the paper. Following Teichm\"uller's notation, we now present the main ideas of this proof. 

Let $\widetilde{\G}^{\prime}$ and $\widetilde{\G}^{\prime\prime}$ be two disjoint simply connected domains containing respectively $0$ and $\infty$ and of reduced modules respectively $\widetilde{M}'$ and $\widetilde{M}''$.  By means of  a homothety of center the origin, one can assume that $z_0 =-1$ does not belong to $\widetilde{\G}^{\prime}\cup \widetilde{\G}^{\prime\prime}$. Therefore, proving the so-called Special Modulastz is equivalent to finding $\delta$ such that  if $\widetilde{M}'+\widetilde{M}^{\prime\prime}\geq -\delta$ then $\widetilde{M}'\leq \varepsilon$ and $\widetilde{M}^{\prime\prime}\leq \varepsilon$.  With this in mind, Teichm\"uller constructs two disjoint simply connected domains that solve an extremal problem. The  construction is as follows.  Let $q>1$ be a parameter. Let $\mathfrak{I}$ and $\mathfrak{A}$ be respectively the interior and exterior of the unit circle in the $w$-plane. By multiplying the argument of a complex number by $q$ one identifies the boundary arc  $w=e^{i\theta}$,  $\left| \theta \right| \leq\frac{\pi}{q}$ of $\mathfrak{I}$ with the boundary curve of $\mathfrak{A}$. Furthermore,  one glues  the remaining boundary arc of $\mathfrak{I}$ by identifying each element with its complex conjugate. Therefore, one obtains a closed Riemann surface of genus $0$ and one can prove (Teichm\"uller does!) that it can be conformally mapped onto the whole $z$-plane where  $w=0$ and $w=\infty$ correspond to respectively $z=0$ and $z=\infty$, the boundary arc $w=e^{i\theta}$,  $\left| \theta \right| \leq\frac{\pi}{q}$ of $\mathfrak{I}$ corresponds to a closed curve  which starts and ends at $z=-q^2$ and encloses the origin, and the remaining boundary arc of $\mathfrak{I}$ corresponds to  the line segment that connects $z=-1$ to $z=-q^2$.  The domains $\mathfrak{I}$ and $\mathfrak{A}$ are then carried to disjoint simply connected domains denoted respectively by $\widetilde{\mathfrak{S}}^{\prime}_q$ and $\widetilde{\mathfrak{S}}^{\prime\prime}_q$. By setting $\widetilde{M}^{\prime}_q$ for the reduced module of $\widetilde{\mathfrak{S}}^{\prime}_q$ and $\widetilde{M}^{\prime\prime}_q$ for the reduced module of $\widetilde{\mathfrak{S}}^{\prime\prime}_q$, Teichm\"uller proves the following extremal property:
\begin{proposition}\label{proposition:extremal-special-modulsatz}
Let $\widetilde{\G}^{\prime}$ and $\widetilde{\G}^{\prime\prime}$ be two disjoint simply connected domains containing respectively $0$ and $\infty$ and both not containing $z=-1$. Let $\widetilde{M}'$ and $\widetilde{M}''$ be their corresponding reduced modules. Then,
$$
\forall q>1, \; q^2 \widetilde{M}' +\widetilde{M}^{\prime\prime}\leq q^2 \widetilde{M}^{\prime}_q+\widetilde{M}^{\prime\prime}_q .
$$
Moreover, the equality holds only when $\widetilde{\G}^{\prime}=\widetilde{\G}^{\prime}_q$ and $\widetilde{\G}^{\prime\prime}=\widetilde{\G}^{\prime}_q$.
\end{proposition}
Teichm\"uller's proof of such a result uses the same idea as the one of his second proof of the so-called Gr\"otzsch module theorem. After showing this result,  Teichm\"uller gives a planar description of  the set $\mathfrak{B}$ of all possible pairs of reduced modules associated with two disjoint simply connected domains defined  as in the statement of Proposition \ref{proposition:extremal-special-modulsatz}. More precisely, he proves that $\mathfrak{B}$ is a convex region that is bounded by the curve formed by the points $\left(\widetilde{M}^{\prime}_q, \widetilde{M}^{\prime\prime}_q \right)$, $\left(0, 0 \right)$ and $\left(\widetilde{M}^{\prime\prime}_q , \widetilde{M}^{\prime}_q\right)$ and when the parameter $q$ varies over the set of real numbers greater than $1$. He continues by determining the exact values of $\widetilde{M}^{\prime}_q$ and $\widetilde{M}^{\prime\prime}_q$ and thus he proves that the boundary curve of $\mathfrak{B}$ is twice differentiable but not thrice differentiable at $(0,0)$.  He  concludes his second proof of Theorem \ref{Specialmodulsatz} by setting 
$$
\begin{cases} \varepsilon = \widetilde{M}^{\prime}_q \\ \delta = - \widetilde{M}^{\prime}_q-\widetilde{M}^{\prime\prime}_q \end{cases}
$$
and he obtains the following estimates:
\begin{equation}\label{deltaepsilon}
\delta\sim \dfrac{\varepsilon ^2}{\log \dfrac{1}{\varepsilon}} , \ \text{as}\ \varepsilon\to 0.
\end{equation}

 \S 4 ends by  discussion of the possibilities of improving $\delta(\varepsilon)$ (making it larger) in the case of the Modulsatz. Teichm\"uller  suggests a method---similar to the one of the proof of the special Modulsatz---for finding the best pair $(M',M'')$ of modules of subdomains $\G'$ and $\G''$ of $\G=\{z \mid r<|z|<R\}$  that would determine such a $\delta$. He does not execute the argument, as it leads to methods which involve elliptic functions, methods that he does not want to use, as he states at several places in his paper.\footnote{A different and shorter proof of a variation of  Theorem \ref{Modulsatz},\index{Modulsatz} together with an estimate equivalent to (\ref{deltaepsilon}),  is due to Pommerenke \cite{Pomm80}, p. 201--202. The result proved by Pommerenke  uses the additional assumption that $\G'\cup\G''=\G$, and thus the estimate (\ref{estimate}) concerns  the points lying on the joint boundary components $\partial \G'\cap \partial \G''$. Pommerenke's proof uses the method of extremal length, properties of univalent conformal mappings, the area theorem, and coefficient estimates. A proof of the Modulsatz, following Pommerenke's method, and further applications can be found in the monograph by  Garnett and Marshall \cite[p. 167--173]{GarMar2005}.}

Teichm\"uller derives in \S 5 a few consequences of the Modulsatz. 
 In particular, he studies the behaviour of a family $\Gamma=\left\lbrace \mathfrak{C}_{\lambda} \right\rbrace_{\lambda \in I}$ of disjoint simple closed curves indexed by a subset $I$ of $\mathbb{R}$ that accumulates at $+\infty$  satisfying the following two properties:
\begin{enumerate}[label=(\roman{*})]
\item\label{item:1} for any pair $(\lambda, \mu)\in I^2$ such that  $\lambda<\mu$,  the closed curve $\mathfrak{C}_{\lambda}$ separates $0$ from $\mathfrak{C}_{\mu}$;
\item\label{item:2} the curve $\mathfrak{C}_{\lambda}$ shrinks to $\infty$ as $\lambda \rightarrow +\infty$.
\end{enumerate}
For such a family of curves, if $\lambda<\mu$ then because of \ref{item:1}. one has a ring domain  bounded by $\mathfrak{C}_{\lambda}$ and $\mathfrak{C}_{\mu}$. If its module is denoted by $M(\lambda, \mu)$ then  \ref{item:2}. implies that, for $\lambda$ fixed, $M(\lambda, \mu)\rightarrow \infty$ as $\mu\rightarrow \infty$.  Teichm\"uller  examines conditions on $M(\lambda,\mu)$ under which the family $\Gamma$ is almost circular at $\infty$, that is, when $\mathfrak{C}_{\mu}$ approaches a circle as $\mu \rightarrow \infty$. Following Teichm\"uller, we recall that $\mathfrak{C}_{\mu}$ approaches a circle as $\mu\rightarrow \infty$ if its so-called \emph{logarithmic oscillation} $\omega(\mu)$ approaches $0$ as $\mu\rightarrow \infty$.  The logarithmic oscillation\index{logarithmic oscillation} $\omega(\mu)$ of $\mathfrak{C}_{\mu}$ is defined as 
$$
 \omega(\mu)=\log \dfrac{r_2(\mu)}{r_1(\mu)}
$$
where $r_1(\mu) =\min\limits_{\C_\mu}|z|$ and $r_2(\mu) =\max\limits_{\C_\mu}|z|$.  The conditions that Teichm\"uller seeks are given by the following result: 
\begin{equation}\label{almost-circularity}
\lim_{\mu\rightarrow \infty}\omega(\mu) =0\iff \lim_{\substack{\chi \rightarrow \infty \\ \chi < \lambda< \mu}}\left\lbrace M(\chi,\mu)-M(\chi,\lambda)-M(\lambda,\mu)\right\rbrace .
\end{equation}  
By observing that for any $\lambda<\mu$, the doubly-connected domain bounded by $\mathfrak{C}_{\lambda}$ and $\mathfrak{C}_{\mu}$ lies between the annuli  $\left\lbrace z \mid r_{2}(\lambda) <\left| z \right| < r_{1}(\mu) \right\rbrace$ and $\left\lbrace z \mid r_{1}(\lambda) <\left| z \right| < r_{2}(\mu) \right\rbrace$ and using (\ref{almost-circularity}), Teichm\"uller proves the sufficient condition for almost circularity. His proof of the necessary condition is more complicated and is based on the Modulsatz,  the normal family argument he obtained in his \S 2.5 and the Koebe distortion theorem.   In addition to finding conditions for the almost circularity of $\Gamma$ at $\infty$,  Teichm\"uller is interested in  the asymptotic behaviour of $r_i (\mu)$ (i=1,2) as $\mu\rightarrow \infty$.  By means of (\ref{almost-circularity}) and properties of modules of ring domains he proves the following:
\begin{lemma}\label{T2} Let
\[
|M(\lambda,\mu)-(\mu-\lambda)|\leq \varphi(\lambda),  \textrm{ for } \lambda<\kappa, \ M(\lambda,\mu)\leq K
\]
 where $ \lim\limits_{\lambda\to \infty} \varphi(\lambda)=0.$ Then for some constant $\alpha$
\[ \lim\limits_{\lambda\to \infty} (\log r_1(\lambda)-\lambda)=\lim\limits_{\lambda\to \infty} (\log r_2(\lambda)-\lambda)=\alpha.\]
 \end{lemma}
As we shall see in the next section, Teichm\"uller uses Lemma \ref{T2} in order to prove the Main Lemma.\index{Main Lemma!Teichm\"uller}\index{Teichm\"uller Main Lemma} Moreover,  as he emphasizes, these investigations on the asymptotic behaviour of the functions $\omega$, $r_1$ and $r_2$ are useful in value distribution theory since they allow to estimate the number of solutions in a disc $\vert z\vert <r$  of an equation of the type $f(z)=a$ where  $f$ is a meromorphic function.  Teichm\"uller  continues by proving in  \S 5 a result similar to Lemma \ref{T2} that  gives  estimates at infinity of the functions $r_{i} \, (i=1,2)$. His proof is mainly based on  the relation (\ref{deltaepsilon}) and, as he points out, the estimates so obtained are actually an improvement of estimates obtained by Ahlfors in \cite[p. 402]{Ahlfors-Acta}. Let us mention that Ahlfors used such estimates to solve a particular case of   Nevanlinna's inverse problem. See \cite{P-value} for further details.

\section{The Main Lemma}\label{s:main}

As mentionned earlier, \S 6 of Teichm\"uller's paper is devoted to the proof of the Main Lemma.\index{Main Lemma!Teichm\"uller}\index{Teichm\"uller Main Lemma} Teichm\"uller  starts this section by defining a quasiconformal mapping\index{quasiconformal mapping}  as a one-to-one continuous mapping $w=w(z)$ between two domains of the complex plane which, except at isolated exceptional points, is continuously differentiable with a nonzero Jacobian (no assumption about orientation-preservation is made). He then recalls that for such a mapping $w=w(z)$ one can define  the \emph{dilatation quotient}\index{dilatation quotient}  as 
\begin{equation}\label{eq:dilatation}
D(z)=D=\dfrac{\Max\left|\dfrac{dw}{dz}\right|}{\Min\left|\dfrac{dw}{dz}\right|},
\end{equation}
where $\dfrac{dw}{dz}$ represents the directional derivative $\lim_{r\rightarrow 0}\frac{w\left(z+re^{i\alpha}\right)-w(z)}{re^{i\alpha}}$ at a point $z$ in the direction $\alpha$, and where the maximum and the minimum are taken over all directions $\alpha$. For a conformal mapping, there is a well-defined notion of derivative which does not depend on the direction. Thus, the directional derivative at each point is constant and the dilatation quotient is equal to $1$. The converse is also true:  at a regular point, if the dilatation quotient is equal to 1, then the derivative does not depend on the direction, which implies that the map is conformal. As recalled by Teichm\"uller, the value $D(z)$ is also equal to the ratio of the major to the minor axis of an ellipse into which an infinitesimal circle is mapped.

Unlike in Teichm\"uller's previous paper \cite{T9}, the dilatation quotient is not assumed to be bounded. He declares that with such an assumption ``one does not necessarily get to the intended function-theoretic applications.'' Indeed, it is now known that this unboundedness condition on the dilatation quotient leads to applications in uniformisation (or the type problem) or in value distribution theory.\index{value distribution} See for instance \cite{GolOst} and also the papers \cite{ABP-type} and \cite{P-value} in the present volume. Let us point out that in the modern definition of quasiconformality, the mapping is assumed to be sense-preserving, the dilatation quotient is assumed to be bounded and the differentiablilty assumption made by Teichm\"uller (with isolated singularities) is replaced by a weaker  assumption of absolute continuity on lines. 

Before Teichm\"uller, Lavrentieff, in his papers  \cite{Lavrentiev1935a, Lavrentiev1935}, already considered such mappings under the term ``almost analytic functions" (cf. the commentary \cite{AP-Lavrentieff} in the present volume).

 One can also add that the study of mappings with unbounded dilatation quotient began to attract significant attention after the publication of the work of David \cite{Dav88}. See also the exposition \cite{Otal} by Otal in Volume III of the present Handbook.  One of the results in David's paper are sufficient  conditions on the dilatation that allow one to solve a  so-called degenerate Beltrami equation. Another approach using estimates of modules of ring domains appears in \cite{Lehto1970,BJ98,GMSV05,BJ04}.   The investigation of the solutions of the degenerate Beltrami equation led to the development of the theory of mappings of exponentially integrable distortion which are not necessarily quasiconformal\index{quasiconformal mapping} in the modern sense. Many of the known properties of  quasiconformal mappings\index{quasiconformal mapping} were extended to this larger class  \cite{AstIwMAr09, IwMAr01}.
After introducing this notion of quasiconformality and before proving the Main Lemma, Teichm\"uller provides an estimate on the distortion of a ring domain under a quasiconformal mapping. More precisely, he proves by means of the length-area method the following result. 
\begin{lemma}\label{T3} Let $r_1<|z|<r_2$ be a circular ring mapped quasiconformally onto a circular ring $\rho_1<|w|<\rho_2$ by a map whose dilatation quotient satisfies
\[ D(z)\leq C(|z|). \]
Then 
\begin{equation}\label{upperlower}
\int\limits_{r_1}^{r_2}\dfrac1{C(r)}\dfrac{dr}{r}\leq \log{\rho_2}-\log {\rho_1}\leq   \int\limits_{r_1}^{r_2}{C(r)}\dfrac{dr}{r}.
\end{equation}
\end{lemma}
 It is worth noting that this result implies a well-known distortion inequality for quasiconformal (in the modern sense) mappings which says that the module of the image of a ring domain by a quasiconformal mapping does not exceed a certain upper bound.
 
As  consequences of such a lemma, Teichm\"uller first  gives a necessary condition on the dilatation quotient for having a quasiconformal mapping that  maps either the punctured plane  onto the unit disc  or  the unit disc onto the punctured plane. More specifically, he  applies (\ref{upperlower}) to show that if there exists a quasiconformal mapping from the punctured plane  onto the unit disc  such that $D\leq C\left( \left| z\right|\right)$, then 
\begin{equation}\label{eq:typecondition}
\int^{\infty}{\frac{dr}{r\cdot C(r)}} <\infty .
\end{equation}
On the other hand, if it is the unit disc that is mapped quasiconformally onto the punctured plane such that $D\leq C\left( \left| z\right|\right)$, then by the same techniques he gets
\begin{equation}\label{eq:typecondition2}
\int^{1}{C(r)\cdot \frac{dr}{r}} =\infty .
\end{equation}
The condition given in (\ref{eq:typecondition}), as pointed out by Teichm\"uller, was already obtained by Lavrentieff in \cite[Th\'{e}or\`{e}me 1]{Lavrentiev1935a}  (see also \cite[Th\'{e}or\`{e}me 7]{Lavrentiev1935}) and as the latter, Teichm\"uller uses it in his \S 7 to obtain a result on the type problem. See \cite{ABP-type} for further details.

Another consequence of Lemma \ref{T3} is the Main Lemma\index{Main Lemma!Teichm\"uller}\index{Teichm\"uller Main Lemma} that we state again for the convenience of the reader.
\begin{theorem}[Main Lemma] \label{ML} Let the punctured z-plane be mapped one-to-one and quasiconformally onto the punctured $w$-plane.  Let the dilatation quotient satisfy $D(z)\leq C(|z|),$ where $C(r)\to 1,$   as $r\to \infty,$ so fast that 
\[
\int\limits^\infty
(C(r)-1)\dfrac{dr}{r}
\]
converges. Then, by approaching infinity
 \[ |w|\sim \mathrm{const} \cdot  |z|.     \]   
 \end{theorem}    
In order to prove this lemma, Teichm\"uller considers  the family of curves $\C_\lambda$ which are the images  of $|z|=e^\lambda$ by the given quasiconformal mapping. For sufficiently large $\lambda$, the family so defined satisfies conditions \ref{item:1}. and \ref{item:2}. stated at the end of \S 5 of this commentary. He uses Lemma \ref{T3}, and Lemma \ref{T2} by setting $\varphi\left( \lambda\right) = \int_{e^{\lambda}}^{\infty}{\left( C(r)-1\right)\frac{dr}{r}}$ to arrive at the conclusion of the statement.



Let us recall that a  proof of this lemma that uses   Ahlfors' distortion theorem was obtained by Wittich in \cite{Wit48a}.  Another proof  can be found in \cite[p. 345ff]{GolOst}.

Later on,  at the end of his most quoted paper \cite{T20}, Teichm\"uller presents a conjecture which says that under the same hypothesis than Theorem \ref{ML}, one has 
 \begin{equation}\label{behavior}
       w(z)\sim \mathrm{const}\cdot  z \textrm{ as }z\rightarrow \infty.
       \end{equation}
In order to justify such a conjecture, Teichm\"uller says that it ``is supported by the `spiraling value distribution' for certain functions and  the induced order of growth. In particular, this holds for the mapping $w=ze^{i\eta(\left| z\right|)}$ ($\eta(r)$ real).'' Teichm\"uller continues by giving an idea of the proof. Later,  a full proof was obtained by Belinski\u{\i} in \cite{Bel53} and Lehto in \cite{Lehto1960}. See also \cite[V.6.6]{LV73}. This result is now referred to as the \emph{Teichm\"uller--Wittich--Belinski\u{\i} theorem}.\index{Teichm\"uller--Wittich--Belinski\u{\i} theorem}\index{Theorem!Teichm\"uller--Wittich--Belinski\u{\i}} A survey of this result with historical comments was written by Drasin in \cite{Drasin1986}. For further applications  of this theorem one refers  to \cite[V.7]{LV73} and \cite[Ch. 11]{BGMR}. One can also find improvements of this theorem  in \cite{ReichWal,BJ94,GM03,MB09,MB10,Shishikura}. 

The last section of Teichm\"uller's paper concerns the type problem and is discussed in our paper \cite{ABP-type} in the present volume.

\medskip

 \noindent {\emph{Acknowledgements.}}  The first and third authors acknowledge support from the U.S. National Science Foundation grants DMS 1107452, 1107263, 1107367 ``RNMS: GEometric structures And Representation varieties'' (the GEAR Network).

 \printindex
 
 {Vincent Alberge, Fordham University, Department of Mathematics, 441 East Fordham Road, Bronx, NY 10458, USA}

\textit{E-mail address}: {\tt valberge@fordham.edu}

 {Melkana Brakalova, Fordham University, Department of Mathematics, 441 East Fordham Road, Bronx, NY 10458, USA}

\textit{E-mail address}: {\tt brakalova@fordham.edu}

 {Athanase Papadopoulos, Institut de Recherche Mathématique Avancée, Universit\'{e} de Strasbourg et CNRS, 7 rue Ren\'{e} Descartes, 67084 Strasbourg Cedex, France}

\textit{E-mail address}: {\tt papadop@math.unistra.fr}


\begin{thebibliography}{1}
 
 \bibitem{abikoff} W. Abikoff, Oswald Teichm\"uller. \emph{The Mathematical Intelligencer} 8 (1986), 8--16.


\bibitem{Ahlfors-CRAS} L. V. Ahlfors, Sur le nombre des valeurs asymptotiques d'une fonction enti\`{e}res d'ordre fini. \emph{C. R. Acad. Sci. (Paris)} 188 (1929), 688--689.

\bibitem{Ahlfors-Lindelof} L. V. Ahlfors, \"{U}ber die asymptotischen Werte der ganzen Funktionen endlicher Ordnung. \emph{Ann. Acad. Sci. Fenn. Ser. A I Math.} 32 (1929), 1--15

\bibitem{Ahlfors:thesis} L. V. Ahlfors, 
Untersuchungen zur Theorie der konformen Abbildung und der ganzen Funktionen.  \emph{Acta Soc. Sci. Fenn.,  Nove ser. A 1} 9 (1930), 1--40.  

 \bibitem{Ahlfors-Acta} L. V. Ahlfors,  \"Uber eine in der neueren Wertverteilungslehre betrachtete Klasse transzendenter Funktionen. \emph{Acta Math.} 58 (1932), 375--406.



 \bibitem{Ahlfors1964} L. V. Ahlfors,  Quasiconformal mappings and their applications. In \emph{Lect. on Modern Math.} (T. L. Saaty, ed.), Volume 2, Wiley, New York 1964, 151--164. 
 
 \bibitem{Ahlfors-QC} L. V. Ahlfors, \emph{Lectures on quasiconformal mappings}. Second Version, with additional chapters by C. J. Earle and I. Kra, M. Shishikura and J. H. Hubbard. American Mathematical Society, Providence 2006.
 
 \bibitem{Ahlfors-CI} L. V. Ahlfors,  \emph{Conformal invariants:  Topics in geometric function theory}. Reprint of the 1973 original. With a foreword by P. Duren, F. W. Gehring and B. Osgood. AMS Chelsea Publishing, Providence, RI, 2010.
 
 \bibitem{Ahlfors1946} L. V. Ahlfors and A. Beurling, Invariants conformes et probl\`emes extr\'emaux. In \emph{Comptes rendus du $10^{me}$ congr\`es des math\'ematiciens scandinaves}, Copenhagen, 1946, 341--351. 
 
 


 \bibitem{Ahlfors-Collected} L. V. Ahlfors, Collected works, in 2 volumes, 
Series: Contemporary Mathematicians Ser., Birk\"auser  Verlag, 1982.



 
 
\bibitem{AP-Lavrentieff}  V. Alberge and A. Papadopoulos, A commentary on Lavrentieff's paper 
  \emph{Sur une classe de repr\'esentations continues}, In \emph{Handbook of Teichm\"uller theory}  (A. Papadopoulos, ed.), Volume VII, EMS Publishing House,  Z\"urich, 2019, p. ???
  


\bibitem{AP1} V. Alberge and A. Papadopoulos, On five papers by Herbert Gr\"otzsch, In \emph{Handbook of Teichm\"uller theory}  (A. Papadopoulos, ed.), Volume VII, EMS Publishing House,  Z\"urich, 2019, p. ???

\bibitem{ABP-type} V. Alberge, M. Brakalova and A. Papadopoulos, Teichm\"uller's work on the type problem, In \emph{Handbook of Teichm\"uller theory}  (A. Papadopoulos, ed.), Volume VII, EMS Publishing House,  Z\"urich, 2019, p. ???

 
  \bibitem{AstIwMAr09}  K. Astala, T. Iwaniec and G. Martin, \emph{Elliptic partial differential equations and quasiconformal mappings in the plane}. Princeton University Press, Princeton 2009.
  
 \bibitem{BGMR} B. Bojarski, V. Gutlyanski\u{\i}, O. Martio, and V. Ryazanov, \emph{Infinitesimal geometry of quasiconformal and bi-Lipschitz mappings in the plane}. EMS tracts in Mathematics 19, EMS Publishing House,  Z\"urich 2013.

\bibitem{Bel53}   P. P. Belinski\u{\i}, Behavior of a quasiconformal mapping at an isolated point (in Russian). \emph{Doklady Akad. Nauk USSR (N.S.)} 91 (1953), 709--710.




\bibitem{Bes52}  A. S. Besicovitch, On two problems of L\"owner. \emph{J. London Math. Soc.} 127 (1952), 141--144.

\bibitem{Beurling-thesis} A. Beurling, \emph{\'{E}tude sur un probl\`{e}me de majoration}. Th\`{e}se pour le doctorat,   Upsal 1933.


 \bibitem{BJ94}  M. A. Brakalova and J. A. Jenkins, On the local behavior of certain homeomorphisms. \emph{Kodai Math. J.} 17 (1994), 201--213.

\bibitem{BJ98} M. A. Brakalova and J. A. Jenkins, On solutions of the Beltrami equation. \emph{J. Anal. Math.} 76 (1998), 67--92.

\bibitem{BJ04} M. A. Brakalova and J. A. Jenkins, On solutions of the Beltrami equation. II. \emph{Publ. Inst. Math. (Beograd) (N.S.)} 75 (2004), 3--8.
 
 \bibitem{MB09}  M. A. Brakalova, \emph{Sufficient and necessary conditions for conformality. Part I. Geometric viewpoint}. \emph{Complex Variables and Elliptic equations} 55 (2010), 137--155.

\bibitem{MB10} M. A. Brakalova, Sufficient and Necessary Conditions for Conformality. Part II. Analytic Viewpoint. \emph{Ann. Acad. Sci. Fenn. Ser. A I Math.} 35 (1988), 235--254.

\bibitem{Caz} C. Andreian Cazacu, Foundations of quasiconformal mapping. In \emph{Handbook of complex analysis: geometric function theory} (R. K\"uhnau, ed.), Volume 2, Elsevier, Amsterdam 2005, 687--753.



\bibitem{Dav88} G. David, Solutions de l'equation de Beltrami avec $\left\| \mu \right\|_{\infty}=1$. \emph{Ann. Acad. Sci. Fenn. Ser. A I Math.} 13 (1988), 25--70.
 
\bibitem{Drasin1977} D. Drasin, The inverse problem of the Nevanlinna theory. \emph{Acta Math.} 138 (1977),  83--151. 
 
\bibitem{Drasin1986} D. Drasin, On the Teichm\"uller--Wittich--Belinskii theorem. \emph{Results in Math.} 10 (1986), 54-65.


\bibitem{DGP} D. Drasin, A. A. Goldberg, and P. Poggi-Corradini, Quasiconformal mappings in value-distribution theory. In \emph{Handbook of complex analysis: geometric function theory} (R. K\"uhnau, ed.), Volume 2, Elsevier, Amsterdam 2005, 755--808. 




\bibitem{GarMar2005}  J. Garnett and D. Marshall, \emph{Harmonic Measure}. Cambridge Univ. Press, Cambridge 2005. 

\bibitem{GolOst} A. A. Goldberg and I. V. Ostrowskii, \emph{Value Distribution of Meromorphic Functions}. Translation of Mathematical Monographs, v. 236, AMS 2008.


     
\bibitem{Groetzsch1928} H. Gr\"otzsch, \"Uber einige Extremalprobleme der konformen Abbildung. \emph{Ber. Verhandl. Sächs. Akad. Wiss. Leipzig Math.-Phys.   Kl.} 80 (1928), 367--376. English translation by A. A'Campo-Neuen, On some extremal problems of conformal mappings. In \emph{Handbook of Teichm\"uller theory} (A. Papadopoulos, ed.), Volume VII, EMS Publishing House,  Z\"urich 2019, p. ??? .

 
\bibitem{Groetzsch1928b} H. Gr\"otzsch,  \"Uber  einige Extremalprobleme der konformen Abbildung.  II. 
 \emph{Ber. Verhandl. Sächs. Akad. Wiss. Leipzig Math.-Phys.   Kl.} 80 (1928), 497--502. English translation by M. Brakalova, On some extremal problems of conformal mappings. II. In \emph{Handbook of Teichm\"uller theory} (A. Papadopoulos, ed.), Volume VII, EMS Publishing House,  Z\"urich 2019, p. ??? .

 
  \bibitem{GM03} V. Gutlyanski\u{\i} and O. Martio, Conformality of a quasiconformal mapping at a point. \emph{J. Anal. Math.} 91 (2003), 179--192.
  
  \bibitem{GMSV05} V. Gutlyanski\u{\i}, O. Martio, T. Sugawa, and M. Vuorinen, On the degenerate Beltrami equation. \emph{Trans. Amer. Math. Soc.} 357 (2005), 875--900.
 
 \bibitem{herschthesis} J. Hersch, Longueurs extr\'{e}males et th\'{e}ories des fonctions. \emph{Comment. Math. Helv.} 29 (1955), 301--337.
 
 \bibitem{hersch55} J. Hersch, Longueurs extr\'{e}males et g\'{e}om\'{e}trie globale. \emph{Annales scientifiques de l'\'{E}.N.S. 3\textsuperscript{e} s\'{e}rie} 72 (1955), 401--414.
 
 \bibitem{HC} A. Hurwitz and R. Courant, \emph{Funktionentheorie}. Springer-Verlag, Berlin 1922.
 
\bibitem{JJ58} J. A. Jenkins, \emph{Univalent Functions and Conformal Mappings}. Ergeb. Math., Heft 18,  Springer Verlag, Berlin--G\"{o}ttingen--Heidelberg 1958.

\bibitem{IwMAr01} T. Iwaniec and G. J. Martin, \emph{Geometric function theory and nonlinear analysis}. Oxford University Press, New York 2001.

\bibitem{JO71} J. A. Jenkins and K. Oikawa, On results of Ahlfors and Hayman. \emph{Illinois J. Math.} 15 (1971), 664--671.

\bibitem{kunzi} H. P. K\"{u}nzi, \emph{Quasikonforme Abbildungen}. Ergeb. Math., Heft 26, Springer-Verlag, Berlin--G\"{o}ttingen--Heidelberg 1960.



   
   
\bibitem{Lavrentiev1935a} M. A. Lavrentieff, Sur une classe de repr\'esentations continues. \emph{C. R. Acad. Sci. (Paris)} 200 (1935), 1010--1013.


\bibitem{Lavrentiev1935} M. A. Lavrentieff, Sur une classe de repr\'esentations continues. \emph{Mat. Sb.} 42 (1935), 407--423. English translation by V. Alberge and A. Papadopoulos, On a class of continuous representations. In \emph{Handbook of Teichm\"uller theory} (A. Papadopoulos, ed.), Volume VII, EMS Publishing House, Z\"urich, to appear. 


\bibitem{Le} T. L\^e V\u{a}n, \"Uber das Umkehrproblem der Wertverteilungstheorie. \emph{Comment. Math. Helv.}  23 (1949), 26--49.

\bibitem{Lethesis} T. L\^e V\u{a}n, Sur un probl\`{e}me d'inversion dans la th\'{e}orie des fonctions m\'{e}romorphes. \emph{Ann. Sci. École Norm. Sup. (3)}  67 (1950), 51--98.

\bibitem{Lehto1960} O. Lehto, On the differentiability of quasiconformal mappings with prescribed complex dilatation.
\emph{Ann. Acad. Sci. Fenn. A I} 275 (1960), 1-28.

\bibitem{Lehto1970} O. Lehto, Homeomorphisms with a given dilatation. \emph{Proceedings of the 15th Scandanavian Congress} 118 (1970), 58--73.







\bibitem{LV73} O. Lehto and L. Virtanen, \emph{Quasiconformal 
mappings in the plane}.  Springer-Verlag, Berlin 1973.



 


\bibitem{N-analytic}  R. Nevanlinna, Analytic functions, Grundlehren der mathematischen Wissenschaften, vol. 162, Springer-Verlag, New York, 1970. Revised translation of \emph{Eindeutige analytische Funktionen}, 2nd ed., Grundlehren der mathematischen Wissenschaften, vol. 46, 1953.




\bibitem{Otal} J. P. Otal, Quasiconformal and BMO-quasiconformal
homeomorphisms, In \emph{Handbook of Teichm\"uller theory}  (A. Papadopoulos, ed.), Volume III, EMS Publishing House,  Z\"urich 2012, 37--70.





\bibitem{P-value} A. Papadopoulos, Value distribution theory and Teichm\"uller's paper \emph{Einfache Beispiele zur Wertverteilungslehre}. In \emph{Handbook of Teichm\"uller theory} (A. Papadopoulos, ed.), Volume VII, EMS Publishing House, Zürich, p. ???

\bibitem{Pomm80} Ch. Pommerenke, Boundary behavior of conformal mappings. In \emph{Aspects of contemporary complex
analysis} (D. A., Brannan and J. Clunie, eds.), Academic Press, New York 1980, 313--332.




\bibitem{ReichWal} E. Reich and H. Walczak, On the Behavior of Quasiconformal Mappings at a Point. \emph{Trans. Amer. Math. Soc.} 117 (1965),  338--351.
 
\bibitem{schappacher&scholtz} N. Schappacher and E. Scholz, Oswald Teichm\"uller -- Leben und Werk. \emph{Jber. d. Dt. Math.-Verein.} 94 (1992), 1--39.

\bibitem{Schottky} F. H. Schottky, \"Uber konforme Abbildung von mehrfach zusammenhängenden Fläche. \emph{J. Reine angew. Math.} 83 (1877), 300--351.

\bibitem{Shishikura} M. Shishikura, Conformality of quasiconformal mappings at a point, revisited. \emph{Ann. Acad. Sci. Fenn. Math.} 43 (2018), 981--990.



\bibitem{T1} O.  Teichm\"uller, Operatoren im Wachsschen Raum, \emph{J. reine angew. Math.} 74 (1935), 73-124.  In \emph{Gesammelte Abhandlungen} (L. V. Ahlfors and F. W. Gehring, eds.), Springer-Verlag, Berlin--Heidelberg--New York 1982, 1--52.
 

  
\bibitem{T8} O. Teichm\"uller, Eine Umkehrung des zweiten Hauptsatzes der Wertverteilungslehre. \emph{Deutsche Math.} 2 (1937), 96--107. In \emph{Gesammelte Abhandlungen} (L. V. Ahlfors and F. W. Gehring, eds.), Springer-Verlag, Berlin--Heidelberg--New York 1982, 158--169.
 

\bibitem{T9} O. Teichm\"uller,  Eine Anwendung quasikonformer Abbildungen auf das Typenproblem.  \emph{Deutsche Math.}  2 (1937), 321--327. In \emph{Gesammelte Abhandlungen} (L. V. Ahlfors and F. W. Gehring, eds.), Springer-Verlag, Berlin--Heidelberg--New York 1982, 171--178. English translation by M. Brakalova, An application of quasiconformal mappings to the type problem. In \emph{Handbook of Teichm\"uller theory} (A. Papadopoulos, ed.), Volume VII, EMS Publishing House, Zürich,  ???

\bibitem{T13} O. Teichm\"uller, Untersuchungen \"uber konforme und quasikonforme Abbildung. \emph{Deutsche Math.} 3 (1938), 621--678. In \emph{Gesammelte Abhandlungen} (L. V. Ahlfors and F. W. Gehring, eds.), Springer-Verlag, Berlin--Heidelberg--New York 1982, 205--262. English translation by M. Brakalova and M. Weiss, Investigations on conformal and quasiconformal mappings. In \emph{Handbook of Teichm\"uller theory} (A. Papadopoulos, ed.), Volume VII, EMS Publishing House, Zürich,  ???


\bibitem{T15} O. Teichm\"uller,  Eine Versch\"arfung des Dreikreisesatzes,  \emph{Deutsche Math.}  4 (1939), 16-22.  In \emph{Gesammelte Abhandlungen} (L. V. Ahlfors and F. W. Gehring, eds.), Springer-Verlag, Berlin--Heidelberg--New York 1982, 276--282.


\bibitem{T17} O. Teichm\"uller, Vermutungen und S\"atze \"uber die Werverteilung gebrochener Funktionen endlicher Ordnung. \emph{Deutsche Math.} 4 (1939), 163-190.  In \emph{Gesammelte Abhandlungen} (L. V. Ahlfors and F. W. Gehring, eds.), Springer-Verlag, Berlin--Heidelberg--New York 1982, 287--314.

\bibitem{T20} O. Teichm\"uller, Extremale quasikonforme Abbildungen und quadratische Differentiale. \emph{Abh. Preu\ss. Akad. Wiss., math.-naturw. K1.} 22 (1939), 1--197.  In \emph{Gesammelte Abhandlungen} (L. V. Ahlfors and F. W. Gehring, eds.), Springer-Verlag, Berlin--Heidelberg--New York 1982, 337--531. English translation by G. Th\'{e}ret, Extremal quasiconformal mappings and quadratic differentials. In \emph{Handbook of Teichm\"uller theory} (A. Papadopoulos, ed.), Volume V, EMS Publishing House, Zürich 2016, 321--483.
  
  
  





\bibitem{T31} O.  Teichm\"uller, 
Ein Verschiebungssatz der quasikonformen Abbildung. \emph{Deutsche Math.} 7 (1944), 336--343.  In \emph{Gesammelte Abhandlungen} (L. V. Ahlfors and F. W. Gehring, eds.), Springer-Verlag, Berlin--Heidelberg--New York 1982, 704--711. English translation  by M. Karbe, A displacement theorem for quasiconformal mapping. In \emph{Handbook of Teichm\"uller theory}  (A. Papadopoulos, ed.), Volume VI, EMS Publishing House,  Z\"urich 2016, 601--608.



%



\bibitem{T33} O. Teichm\"uller, Einfache Beispiele zur Wertverteilungslehre. \emph{Deutsche Math.}  7 (1944), 360--368. In \emph{Gesammelte Abhandlungen} (L. V. Ahlfors and F. W. Gehring, eds.), Springer-Verlag, Berlin--Heidelberg--New York 1982, 728--736. English translation by A. A'Campo-Neuen, Simple examples in value distribution theory.  In \emph{Handbook of Teichm\"uller theory} (A. Papadopoulos, ed.), Volume VII, EMS Publishing House, Zürich, p. ???







\bibitem{Ullrich} E. Ullrich, Zum Umkehrproblem der Wertvertelungslehre. \emph{Nachr. Ges. Wiss G\"ottingen, Math.-Phys. Kl. I, N. F.}, 1 (1936), 135--150.


\bibitem{Wit48a} H. Wittich, Zum Beweis eines Satzes \"uber quasikonforme Abbildungen. \emph{Math. Z.} 51,
(1948), 278--288 .

\bibitem{Wit68} H. Wittich, \emph{Neuere Untersuchungen \"uber eindeutige analytische Funktionen}. Ergeb. Math., Heft 8, Springer-Verlag, Berlin--G\"{o}ttingen--Heidelberg 1955.
\end{thebibliography}
\end{document}